\documentclass{article}
\usepackage{geometry}
\usepackage{amsmath}
\usepackage{amssymb}
\usepackage{bm}
\usepackage{graphicx}
\graphicspath{{./figures/}}
\usepackage{color}
\usepackage{tabu}
\usepackage{booktabs}
\usepackage{authblk}
\usepackage{subfigure}
\usepackage{enumerate}
\usepackage{tabu}
\usepackage{multirow}
\usepackage{diagbox}
\usepackage{amsopn}
\usepackage{amsfonts}
\usepackage{epstopdf}
\usepackage{algorithmic}
\usepackage{epstopdf}

\newtheorem{remark}{Remark}

\newcommand{\sS}{\mathbb{S}}
\newcommand{\eps}{\varepsilon}

\newcommand*\diff{\mathop{}\!\mathrm{d}}

\begin{document}

\title{Asymptotic-Preserving Neural Networks based on Even-odd Decomposition for Multiscale Gray Radiative Transfer Equations}
\author[1,2]{Keke Wu}
\author[3]{Xizhe Xie}
\author[3]{Wengu Chen}
\author[3]{Han Wang}
\author[4,5, \thanks{Corresponding author: zhengma@sjtu.edu.cn}]{Zheng Ma}

\affil[1]{School of Mathematical Sciences and Suzhou Institute for Advanced Research, University of Science and Technology of China, Jiangsu,
    215217, P. R. China.}
\affil[2]{Suzhou Big Data \& AI Research and Engineering Center}
\affil[3]{Institute of Applied Physics and Computational Mathematics, Beijing, 100088,  P. R. China.}
\affil[4]{School of Mathematical Sciences, Shanghai Jiao Tong University, Shanghai,
200240, P. R. China.}
\affil[5]{Institute of Natural Sciences, MOE-LSC, Shanghai Jiao Tong University, Shanghai, 200240, P. R. China.}

\date{\today}

\maketitle

\begin{abstract}
We present a novel Asymptotic-Preserving Neural Network (APNN) approach utilizing even-odd decomposition to tackle the nonlinear gray radiative transfer equations (GRTEs). Our AP loss demonstrates consistent stability concerning the small Knudsen number, ensuring the neural network solution uniformly converges to {the diffusion limit solution}. This APNN method alleviates the rigorous conservation requirements while simultaneously incorporating an auxiliary deep neural network, distinguishing it from the APNN method based on micro-macro decomposition for GRTE. Several numerical problems are examined to demonstrate the effectiveness of our proposed APNN technique.
\end{abstract}



\section{Introduction}
Deep learning methods and deep neural networks (DNNs) have attracted considerable attention in the scientific community, particularly in the context of resolving partial differential equations (PDEs)~\cite{beck2020,cai2021least,E2018,liao2019deep,raissi2019physics,deepGalerkin2018,zang2020weak}. 
The primary motivation behind these approaches is to represent the solutions of PDE problems using deep neural networks. 
This results in high-dimensional and nonconvex minimization problems, setting them apart from classical numerical methods. 
A significant advantage of deep learning methods is their mesh-free nature, allowing them to handle PDEs in complex domains and geometries. 
Moreover, they offer flexibility and ease of implementation, making them suitable for tackling high-dimensional problems.
Despite these advantages, deep learning methods do come with some potential drawbacks, including lengthy training times, convergence challenges, and reduced accuracy. However, the concept of operator learning provides a promising solution for a class of PDEs by training the neural network just once~\cite{li2020fourier,lu2021learning,zhang2021mod,li2021physics,wang2021learning,xiong2023koopman,xu2023transfer,cao2023lno,liu2022ht,wu2023asymptotic}. 
It is worth noting that certain aspects related to the convergence theory of these methods still require further clarification.
To explore alternative machine learning approaches for solving partial differential equations, we recommend interested readers to consult the exemplary review article~\cite{beck2020}.

In recent years, there has been extensive research on employing deep neural networks to tackle multiscale kinetic equations and hyperbolic systems~\cite{CLM,HJJL,lee2021model,wuAPNN,lu2022solving,li2022model,bertaglia2022asymptotic1,bertaglia2022asymptotic2,lou2021physics,wu2023asymptotic,wuAPNNv2,zhang2023asymptotic,lee2024structure,li2024solving,li2024macroscopic,liu2025}. 
These problems, characterized by features at multiple scales, have gained significant importance in various scientific investigations.
When dealing with partial differential equations (PDEs), there are several choices available to formulate the loss function, such as variational formulation (DRM), least-squares formulation (PINN, DGM), weak formulation (WAN), and others. 
However, when addressing multiscale kinetic equations, the conventional Physics-Informed Neural Networks (PINNs) may encounter instability due to the presence of small scales~\cite{lu2022solving,li2022model,wuAPNN,wuAPNNv2,wu2023asymptotic}.
A crucial consideration is determining what constitutes a ``good'' loss function, one that accounts for conservation, symmetry, parity, and other essential factors. 
An approach to resolve multiscale kinetic equations using deep neural networks is to design a loss that captures the limiting macroscopic behavior, referred to as Asymptotic-Preserving (AP) loss, justifying the use of Asymptotic-Preserving Neural Networks (APNNs)~\cite{wuAPNN}. 
Accordingly, we proposed an APNN method for time-dependent linear transport equations with diffusive scaling and uncertainties, based on micro-macro decomposition, demonstrating that the loss exhibits AP behavior with respect to the Knudsen number approaching zero.
It is crucial to note that the APNN method, based on micro-macro decomposition, imposes strict conservation prerequisites. 
Failure to meet these criteria may lead to imprecise outcomes in the deep neural network approximation. Consequently, we have made progress towards refining the APNN approach by utilizing an even-odd decomposition to address time-dependent linear transport equations, thereby relaxing the stringent conservation prerequisites~\cite{wuAPNNv2}.

The gray radiative transfer equations (GRTEs) play a crucial role in modeling photon transportation and energy interactions with surrounding substances, finding diverse applications in fields such as astrophysics, inertial or magnetic confinement fusion, high-temperature flow systems, and more~\cite{chandrasekhar2013radiative,siegel2001thermal,sun2015asymptotic}. 
However, achieving precise simulations of GRTEs presents a formidable challenge due to their complex nature characterized by high dimensionality, strong coupling nonlinearity, and multiscale features arising from varying opacities exhibited by background materials.
In a recent study~\cite{li2022model}, the authors proposed a model-data asymptotic-preserving neural network method based on micro-macro decomposition for GRTEs, demonstrating promising results in simulating the nonlinear non-stationary behavior of these equations. 
Building upon this valuable work, we have introduced a novel APNN method based on even-odd decomposition to solve the nonlinear GRTEs. 
Our approach features a novel AP loss that exhibits uniform stability concerning the small Knudsen number, ensuring that the neural network solution converges uniformly to {the diffusion limit solution}. 
Notably, this APNN method relaxes the stringent conservation prerequisites while simultaneously introducing an auxiliary deep neural network, setting it apart from the APNN method based on micro-macro decomposition for GRTEs.

This paper comprises the following outline:
In Section 2, a detailed illustration of Asymptotic-Preserving Neural Networks based on micro-macro decomposition and even-odd decomposition for gray radiative transfer equations and the construction of the AP loss functions are given.
Numerous numerical examples are presented in Section 3 to demonstrate the effectiveness for both APNN methods.
The summary of the contributions and findings of the paper is concluded in Section 4.

\section{Methodology}

\subsection{The gray radiative transfer equation and its diffusion limit }
Consider the scaled form of the gray radiative transfer equations in a bounded domain $\tau \times D \times \sS^2$~\cite{li2022model}:
\begin{equation}\label{eq:GRTE}
    \left\{
    \begin{aligned}
         & \frac{\eps^2}{c} \frac{\partial I}{\partial t}+\eps \Omega \cdot \nabla I =\sigma \left(\frac{1}{4 \pi}acT^4-I\right),    \\
         & \eps^2 C_v\frac{\partial T}{\partial t} = \sigma \left(\int_{\sS^2} I \diff{\Omega} -acT^4\right), \\
         & \mathcal{B} I = 0, \; \\
         & I(t=0, x, \Omega) = I_0(x, \Omega), \\
         & T(t=0, x) = T_0(x),
    \end{aligned}
    \right.
\end{equation}
where $I(t,x,\Omega)$ is the radiation intensity at time $t \in \tau$, space point $x \in D$, and angular direction $\Omega \in \sS^2$, $T(t,x)$  
is the material temperature, $a, c, \sigma$ denote the radiation constant, scaled speed of light, opacity and $C_v$ is the scaled heat capacity. $ \mathcal{B}$ is the boundary operator for $I$. The parameter $\eps > 0$ is called the Knudsen number which characterizes the ratio of mean free path over the system characteristic length. 

In the context of a one-dimensional scenario, the GRTE simplifies to
\begin{equation}\label{eq:GRTE1d}
    \left\{
    \begin{aligned}
         & \frac{\eps^2}{c} \frac{\partial I}{\partial t}+\eps \mu \frac{\partial I}{\partial x} =\sigma \left(\frac{1}{2}acT^4 -I\right), \; (t,x,\mu) \in \tau \times D \times [-1,1],               \\
         & \eps^2 C_v\frac{\partial T}{\partial t} = \sigma \left(\int_{-1}^{1}I \diff{\mu} -acT^4 \right), \; (t,x) \in \tau \times D.
    \end{aligned}
    \right.
\end{equation}
When the temperature of the material aligns with the temperature of radiation, expressed as $T_r = (\frac{1}{ac}\int I \diff{\Omega} )^{1/4}$, equation \eqref{eq:GRTE} transitions into the scaled linear transport model
\begin{equation}\label{eq:LTE2d}
\frac{\eps^2}{c}\partial_t I + \eps \Omega \cdot \nabla I = \sigma\left(\frac{1}{4\pi}\int_{S^2} I \diff{\Omega} -I\right),
\end{equation}
{in 1d case}, one can obtain
\begin{equation}\label{eq:LTE1d}
\frac{\eps}{c} \frac{\partial I}{\partial t} + \mu \frac{\partial I}{\partial x} = \frac{\sigma}{\eps} \left(\frac{1}{2}\int_{-1}^{1}I \diff{\mu} -I\right).
\end{equation}
Equation~\eqref{eq:GRTE} represents a relaxation model about the radiation intensity within the context of local thermodynamic equilibrium, with the emission source originating from the background medium, as dictated by the Planck function corresponding to the local material temperature, more precisely, $\sigma acT^4 / {4\pi}$.
As the parameter {$\eps$} tends towards zero, while disregarding boundaries and initial moments, the radiation intensity denoted as $I$ converges towards a Planck function at the local temperature~\cite{larsen1983asymptotic}. This can be stated as $I^{(0)}= ac\left(T^{(0)}\right)^4 / 4 \pi$. Additionally, the local temperature $T^{(0)}$ satisfies a diffusion equation:
\begin{equation}\label{eq:diffusion}
\frac{\partial}{\partial t}\left(C_v T^{(0)}\right) + a \frac{\partial}{\partial t}\left(T^{(0)}\right)^4 = \nabla \cdot \frac{ac}{3\sigma}
\nabla \left(T^{(0)}\right)^4.
\end{equation}

Inspired by the recent works of \cite{lu2022solving, li2022model, wuAPNNv2} for solving gray radiative transfer equations based on micro-macro decomposition, we shall derive the corresponding APNN method for GRTEs based on even-odd decomposition.

\subsection{Micro-macro decomposition method for GRTEs}
First, let's review the APNN method based on micro-macro decomposition to solve GRTEs\cite{lu2022solving, li2022model}.
In this paper, we consider the 1d case for GRTEs under diffusive scaling and given by 
\begin{equation*}
    \left\{
    \begin{aligned}
         & \frac{\eps^2}{c} \frac{\partial I}{\partial t}+\eps \mu \frac{\partial I}{\partial x} =\sigma \left(\frac{1}{2}acT^4 -I\right), \; (t,x,\mu) \in \tau \times D \times [-1,1],               \\
         & \eps^2 C_v\frac{\partial T}{\partial t} = \sigma \left(\int_{-1}^{1}I \diff{\mu} -acT^4 \right), \; (t,x) \in \tau \times D.
    \end{aligned}
    \right.
\end{equation*}

Decompose the radiative intensity $I(t,x,\mu)$ into the equilibrium $\rho(t,x)$ and the non-equilibrium part $g(t,x,\mu)$ as follows:
\begin{equation*}
    I(t,x,\mu) = \rho(t,x) + \frac{\eps}{\sqrt{\sigma_0}}  g(t,x,\mu), \; \left \langle g \right \rangle := \frac{1}{2} \int_{-1}^{1} g(t, x, \mu) \diff{\mu} = 0.    
\end{equation*}
Here, $\sigma_0 > 0$ is a reference opacity constant, {typically taken as the maximum value of opacity~\cite{xiong2022high}.
The inclusion of $\sqrt{\sigma_0}$ is due to the consideration that when $\sigma$ is large, the original equation approaches thermodynamic equilibrium.}

Let us define the operator {$\Pi: \Pi(h)(\mu) = \langle h \rangle = \frac{1}{2} \int_{-1}^{1} h(t, x, \mu) \diff{\mu}$} and the identity operator $\text{Id}$.  
By using these operators, one can derive the micro-macro system for the gray radiative transfer equations:
\begin{equation}\label{eq:GRTE-micro-macro}
    \left\{
    \begin{aligned}
        & \frac{1}{c} \partial_t \rho  + \frac{1}{\sqrt{\sigma_0}} \left \langle  \mu \cdot \partial_x g  \right \rangle =  - \frac{1}{2} C_v \partial_t T, \\     
        & \frac{\eps^2}{c} \partial_t g + \eps \left ( {\text{Id}} - \Pi \right ) \left ( \mu \cdot \partial_x g \right ) + \sqrt{\sigma_0} \mu \cdot \partial_x \rho   =  - \sigma g, \\
        & \eps^2 C_v \partial_t T = \sigma \left ( 2 \rho - acT^4 \right ).
    \end{aligned}
    \right.
\end{equation}

{When $\eps \to 0$, the system formally reduces to its asymptotic limit form, leading to the following equations for the gray radiative transfer equations:}

\begin{equation*}
    \left\{
    \begin{aligned}
        & \frac{1}{c} \partial_t \rho  + \frac{1}{\sqrt{\sigma_0}} \left \langle  \mu \cdot \partial_x g  \right \rangle =  - \frac{1}{2} C_v \partial_t T, \\     
        & \sqrt{\sigma_0} \mu \cdot \partial_x \rho   =  - \sigma g, \\
        & 0 = \sigma \left ( 2 \rho - acT^4 \right ),
    \end{aligned}
    \right.
\end{equation*}
which exactly results in the nonlinear diffusion limit equation
\begin{equation*}
\frac{\partial}{\partial t}\left(C_v T\right) + a \frac{\partial}{\partial t}\left(T\right)^4 = \nabla \cdot \frac{ac}{3\sigma}
\nabla \left(T\right)^4.
\end{equation*}

\subsection{Even-odd decomposition method for GRTEs~{\cite{wuAPNNv2, jin2000uniformly}}}
Next, we shall continue our investigation in the realm of GRTEs under the influence of diffusive scaling, focusing exclusively on the one-dimensional scenario, as represented by the following expression:
\begin{equation*}
    \left\{
    \begin{aligned}
         & \frac{\eps^2}{c} \frac{\partial I}{\partial t}+\eps \mu \frac{\partial I}{\partial x} =\sigma \left(\frac{1}{2}acT^4 -I\right), \; (t,x,\mu) \in \tau \times D \times [-1,1],               \\
         & \eps^2 C_v\frac{\partial T}{\partial t} = \sigma \left(\int_{-1}^{1}I \diff{\mu} -acT^4 \right), \; (t,x) \in \tau \times D.
    \end{aligned}
    \right.
\end{equation*}
By decomposing the equation and establishing distinct even and odd parities, we can proceed as follows~\cite{jin2000uniformly}:
\begin{equation}
    \begin{aligned}
        r(t, x, \mu) & = \frac{1}{2}[I(t, x, \mu) + I(t, x, -\mu)], \; 0 \le \mu \le 1,             \\
        j(t, x, \mu) & = \frac{\sqrt{\sigma_0}}{2\eps}[I(t, x, \mu) - I(t, x, -\mu)],  \; 0 \le \mu \le 1,
    \end{aligned}
\end{equation}
here, we set $\sigma_0 > 0$ as a reference opacity, maintaining its significance as in the previous context and one can obtain
\begin{equation}
    \left\{
    \begin{aligned}
         & \frac{\eps^2}{c} \partial_t r + \frac{\eps^2}{\sqrt{\sigma_0}} \cdot \mu \partial_x j = \sigma \left ( \frac{1}{2} a c T^4 - r \right ), \\
         & \frac{\eps^2}{c \sqrt{\sigma_0}} \partial_t j + \mu \partial_x r = -\frac{\sigma}{\sqrt{\sigma_0}} j.
    \end{aligned}
    \right.
\end{equation}

Let us denote $\rho = \left \langle r \right \rangle = \int_0^1 r(t, x, \mu) \diff{\mu}$. Integrating over $\mu$, the first equation yields:
\begin{equation}
     \frac{\eps^2}{c} \partial_t \rho  + \frac{\eps^2}{\sqrt{\sigma_0}} \cdot  \left \langle  \mu \partial_x j \right \rangle = \sigma \left ( \frac{1}{2} a c T^4 - \rho \right ).
\end{equation}
Finally, the GRTEs can be reformulated as an even-odd system:
\begin{equation}\label{eq:even-odd}
    \left\{
    \begin{aligned}
         & \frac{\eps^2}{c} \partial_t r + \frac{\eps^2}{\sqrt{\sigma_0}} \cdot \mu \partial_x j = \sigma \left ( \frac{1}{2} a c T^4 - r \right ), \\
         & \frac{\eps^2}{c \sqrt{\sigma_0}} \partial_t j + \mu \partial_x r = -\frac{\sigma}{\sqrt{\sigma_0}} j, \\
         &  \frac{1}{c} \partial_t \rho  + \frac{1}{\sqrt{\sigma_0}} \cdot \left \langle  \mu \partial_x j \right \rangle = - \frac{1}{2} C_v\frac{\partial T}{\partial t}, \\
         & \eps^2 C_v\frac{\partial T}{\partial t} = \sigma \left(2 \rho -acT^4 \right), \\
         & \rho = \left \langle r \right \rangle.
    \end{aligned}
    \right.
\end{equation}

As $\eps$ tends to zero, the aforementioned equation formally converges to
\begin{equation}
    \left\{
    \begin{aligned}
         & \frac{1}{2} a c T^4 = r,   \\
         & \mu \partial_x r = -\frac{\sigma}{\sqrt{\sigma_0}} j,   \\
         & \frac{1}{c} \partial_t \rho  + \frac{1}{\sqrt{\sigma_0}} \cdot  \left \langle  \mu \partial_x j \right \rangle = - \frac{1}{2} C_v\frac{\partial T}{\partial t}, \\
         & \rho = \frac{1}{2} acT^4.
    \end{aligned}
    \right.
\end{equation}
Upon substituting these equations into the third equation, we will obtain:
\begin{equation*}
    \frac{\partial}{\partial t}\left(C_v T \right) + \frac{\partial}{\partial t} \left ( a T^4 \right ) = \nabla \cdot \frac{c}{3\sigma}
\nabla \left ( a T^4 \right ),
\end{equation*}
yielding precisely the diffusion equation~\eqref{eq:diffusion}.

\subsection{Solving GRTEs by DNNs}
Firstly, let us introduce the conventional notations commonly used for DNNs
\footnote{BAAI.2020.\ Suggested Notation for Machine Learning.\ https://github.com/mazhengcn/suggested-notation-for-machine-learning.}.
An $L$-block ResNet~\cite{he2016deep} is recursively defined as follows:
\begin{equation}
    \begin{aligned}
        f_{\theta}^{[0]}(x) & = W^{[0]} x + b^{[0]},                                                                                                                                  \\
        f_{\theta}^{[l]}(x) & =  f_{\theta}^{[l-1]}(x) + \sigma \circ( W_2^{[l-1]} \sigma \circ (W_1^{[l-1]} f_{\theta}^{[l-1]}(x) + b_1^{[l-1]}) + b_2^{[l-1]}), \, 1 \le l \le L-1, \\
        f_{\theta}(x)       & = f_{\theta}^{[L]}(x) = W^{[L-1]} f_{\theta}^{[L-1]}(x) + b^{[L-1]},
    \end{aligned}
\end{equation}
here, we have $W_1^{[l]}, W_2^{[l]} \in \mathbb{R}^{m_{l+1}\times m_l}$, and $b_1^{l}, b_2^{l} \in \mathbb{R}^{m_{l+1}}$. The dimensions are defined as follows: $m_0 = d_{in} = d$, representing the input dimension, and $m_L = d_0$, representing the output dimension. Additionally, $\sigma$ denotes a scalar function, and "$\circ$" denotes an entry-wise operation.
In general, ResNet is composed of several residual blocks, with each block consisting of one input, two weight layers, two activation functions, one identical (shortcut) connection, and one output.
We represent the set of parameters as $\theta$.

To solve the gray radiative transfer equation using DNNs, a representative approach is known as Physics-Informed Neural Networks. The vanilla PINN method formulates the PINN loss as the least squares of the residual of the GRTE, in combination with boundary and initial conditions. In the PINN method, it is often necessary to employ two DNNs to parameterize $I(t, x, \mu)$ and $T(t, x)$, respectively:
\begin{equation}
    I^{\text{NN}}_{\theta}(t, x, \mu) := \sigma^+ \left(\tilde{I}^{\text{NN}}_{\theta}(t, x, \mu) \right) \approx I(t, x, \mu),
\end{equation}
\begin{equation}
    T^{\text{NN}}_{\theta}(t, x) := \sigma^+\left(\tilde{T}^{\text{NN}}_{\theta}(t, x)\right) \approx T(t, x),
\end{equation}
which $\sigma^+$ is a scalar function to keep them positive if necessary. 

The PINN empirical loss for the GRTE can be expressed as follows:
\begin{equation}
    \mathcal{R}^{\eps}_{\text{PINN}} = \mathcal{R}^{\eps}_{\text{residual}} + \mathcal{R}^{\eps}_{\text{boundary}} + \mathcal{R}^{\eps}_{\text{initial}},
\end{equation}
where $\mathcal{R}^{\eps}_{\text{residual}}, \mathcal{R}^{\eps}_{\text{boundary}}, \mathcal{R}^{\eps}_{\text{initial}}$ are denoted by
\begin{equation}
    \begin{aligned}
        \mathcal{R}^{\eps}_{\text{residual}}   & = \; \frac{1}{N_\text{int}} \sum_{i=1}^{N_\text{int}} | \frac{\eps^2}{c} \partial_t I^{\text{NN}}_{\theta}(t_i,x_i,\mu_i) + \eps \mu_i \partial_x I^{\text{NN}}_{\theta}(t_i,x_i,\mu_i) \\
        & \quad \quad \quad \quad \quad \quad
        - \sigma ( \frac{1}{2} ac (T^{\text{NN}}_{\theta}(t_i,x_i))^4 - I^{\text{NN}}_{\theta}(t_i,x_i,\mu_i)) |^2                                                                                                                                                       \\
                                                      & + \frac{1}{N_\text{int}} \sum_{i=1}^{N_\text{int}}  |\eps^2 C_v \partial_t T^{\text{NN}}_{\theta}(t_i,x_i) - \sigma (\int_{-1}^1 I^{\text{NN}}_{\theta}(t_i,x_i,\mu) \diff{\mu} - ac (T^{\text{NN}}_{\theta}(t_i,x_i))^4) |^2,        \\
                \mathcal{R}^{\eps}_{\text{boundary}}   & = \frac{1}{N_\text{bdy}} \sum_{i=1}^{N_\text{bdy}}  |\mathcal{B} I^{\text{NN}}_{\theta}(t_i,x_i,\mu_i)|^2, \\
        \mathcal{R}^{\eps}_{\text{initial}}    & = \frac{1}{N_0}\sum_{i=1}^{N_0}  |T^{\text{NN}}_{\theta}(0, x_i) -  T_0 (x_i)|^2 + \frac{1}{N_0} \sum_{i=1}^{N_0} |I^{\text{NN}}_{\theta}(0, x_i, \mu_i) - I_{0}(x_i,\mu_i)|^2. 
    \end{aligned}
\end{equation}
Here, $N_\text{int}, N_\text{bdy}, N_0$ are the number of sample points of corresponding domains.

For the APNN method based on micro-macro decomposition, three DNNs are employed to parameterize $\rho(t, x)$, $T(t, x)$, and $g(t, x, \mu)$. Therefore, three networks are utilized in this approach:
\begin{equation}
    \rho^{\text{NN}}_{\theta}(t, x) := \sigma^+\left( \tilde{\rho}^{\text{NN}}_{\theta}(t, x)\right) \approx \rho(t, x),
\end{equation}
\begin{equation}
    T^{\text{NN}}_{\theta}(t, x) := \sigma^+\left( \tilde{T}^{\text{NN}}_{\theta}(t, x)\right) \approx T(t, x),
\end{equation}
\begin{equation}\label{eq:g-net}
    g^{\text{NN}}_{\theta}(t, x, \mu) := \tilde{g}^{\text{NN}}_{\theta}(t, x, \mu) - \langle \tilde{g}^{\text{NN}}_{\theta} \rangle (t, x) \approx g(t, x, \mu).
\end{equation}
Notice that 
\begin{equation*}
    \langle g^{\text{NN}}_{\theta} \rangle = \langle \tilde{g}^{\text{NN}}_{\theta} \rangle - \langle \tilde{g}^{\text{NN}}_{\theta} \rangle = 0, \quad \forall \; t, x.
\end{equation*}

\begin{remark}
 In Eq. \eqref{eq:g-net}, we highlighted the construction of $g^{\text{NN}}_{\theta}$, which is essential in formulating the APNN loss, as illustrated in \cite{wuAPNN}. When considering $\langle g \rangle = 0$, \cite{li2022model} introduces this condition as a soft constraint in the loss of APNN.
\end{remark}

The APNN empirical risk for the GRTE based on micro-macro decomposition can be represented as:
\begin{equation}
    \mathcal{R}^{\eps}_{\text{APNN-MM}} = \mathcal{R}^{\eps}_{\text{residual}} + \mathcal{R}^{\eps}_{\text{boundary}} + \mathcal{R}^{\eps}_{\text{initial}},
\end{equation}
where $\mathcal{R}^{\eps}_{\text{residual}}, \mathcal{R}^{\eps}_{\text{boundary}}, \mathcal{R}^{\eps}_{\text{initial}}$ are denoted by
\begin{equation}
    \begin{aligned}
        \mathcal{R}^{\eps}_{\text{residual}}   & = \; \frac{1}{N_\text{int}} \sum_{i=1}^{N_\text{int}} |   \frac{1}{c} \partial_t \rho^{\text{NN}}_{\theta}(t_i,x_i)  + \frac{1}{\sqrt{\sigma_0}} \left \langle  \mu \cdot \partial_x g^{\text{NN}}_{\theta}  \right \rangle (t_i,x_i) + \frac{1}{2} C_v \partial_t T^{\text{NN}}_{\theta}(t_i,x_i) |^2                  \\
                                                      & + \frac{1}{N_\text{int}} \sum_{i=1}^{N_\text{int}}  | \frac{\eps^2}{c} \partial_t g^{\text{NN}}_{\theta}(t_i,x_i,\mu_i) + \eps \left ( {\text{Id}} - \Pi \right ) \left ( \mu \cdot \partial_x g^{\text{NN}}_{\theta} \right )(t_i,x_i,\mu_i)  \\
                                                      & \quad \quad \quad \quad \quad \quad + \sqrt{\sigma_0} \mu_i \cdot \partial_x \rho^{\text{NN}}_{\theta}(t_i,x_i) + \sigma g^{\text{NN}}_{\theta}(t_i,x_i,\mu_i) |^2       \\
                                                      & +\frac{1}{N_\text{int}}\sum_{i=1}^{N_\text{int}} | \eps^2 C_v\frac{\partial T^{\text{NN}}_{\theta}}{\partial t}(t_i,x_i) - \sigma (2 \rho^{\text{NN}}_{\theta}(t_i, x_i) - a c (T^{\text{NN}}_{\theta}(t_i,x_i))^4 ) |^2 ,            \\
        \mathcal{R}^{\eps}_{\text{boundary}}   & = \frac{1}{N_\text{bdy}} \sum_{i=1}^{N_\text{bdy}}  |\mathcal{B}(\rho^{\text{NN}}_{\theta} + \frac{\eps}{\sqrt{\sigma_0}} g^{\text{NN}}_{\theta})(t_i,x_i,\mu_i) |^2, \\
                \mathcal{R}^{\eps}_{\text{initial}}    & = \frac{1}{N_0}\sum_{i=1}^{N_0}  | T^{\text{NN}}_{\theta}(0, x_i) - T_0 (x_i)|^2 + \frac{1}{N_0} \sum_{i=1}^{N_0} | (\rho^{\text{NN}}_{\theta} + \frac{\eps}{\sqrt{\sigma_0}} g^{\text{NN}}_{\theta})(0, x_i,\mu_i) - I_{0}(x_i,\mu_i) |^2.
    \end{aligned}
\end{equation}

For the APNN method based on even-odd decomposition, four DNNs are utilized to parameterize $\rho(t, x)$, $T(t, x)$, $r(t, x, \mu)$, and $j(t, x, \mu)$, respectively. Hence, four networks are employed in this approach:
\begin{equation}
    \rho^{\text{NN}}_{\theta}(t, x) := \sigma^+\left( \tilde{\rho}^{\text{NN}}_{\theta}(t, x)\right) \approx \rho(t, x),
\end{equation}
\begin{equation}
    T^{\text{NN}}_{\theta}(t, x) := \sigma^+\left( \tilde{T}^{\text{NN}}_{\theta}(t, x)\right) \approx T(t, x),
\end{equation}
\begin{equation}
    r^{\text{NN}}_{\theta}(t, x, \mu) := \sigma^+ \left(\tilde{r}^{\text{NN}}_{\theta}(t, x, \mu) + \tilde{r}^{\text{NN}}_{\theta}(t, x, -\mu) \right) \approx r(t, x, \mu),
\end{equation}
\begin{equation}
    j^{\text{NN}}_{\theta}(t, x, \mu) :=
    \tilde{j}^{\text{NN}}_{\theta}(t, x, \mu) - \tilde{j}^{\text{NN}}_{\theta}(t, x, -\mu)  \approx j(t, x, \mu).
\end{equation}
Indeed, through the APNN method based on even-odd decomposition, one can automatically ensure that $r^{\text{NN}}_{\theta}(t, x, \mu)$ and $j^{\text{NN}}_{\theta}(t, x, \mu)$ satisfy the even-odd properties without any explicit constraints.
 
The APNN empirical risk for the GRTE based on even-odd decomposition is
\begin{equation}
    \mathcal{R}^{\eps}_{\text{APNN-EO}} = \mathcal{R}^{\eps}_{\text{residual}} + \mathcal{R}^{\eps}_{\text{constraint}} + \mathcal{R}^{\eps}_{\text{boundary}} + \mathcal{R}^{\eps}_{\text{initial}},
\end{equation}
where $\mathcal{R}^{\eps}_{\text{residual}}, \mathcal{R}^{\eps}_{\text{constraint}}, \mathcal{R}^{\eps}_{\text{boundary}}, \mathcal{R}^{\eps}_{\text{initial}}$ are denoted by
\begin{equation}
    \begin{aligned}
        \mathcal{R}^{\eps}_{\text{residual}}   & = \; \frac{1}{N_\text{int}} \sum_{i=1}^{N_\text{int}} | \frac{\eps^2}{c} \partial_t r^{\text{NN}}_{\theta}(t_i,x_i,\mu_i) + \frac{\eps^2}{\sqrt{\sigma_0}} \mu_i \partial_x j^{\text{NN}}_{\theta}(t_i,x_i,\mu_i)                                                                                     \\
                                                      & \quad \quad \quad \quad \quad \quad - \sigma ( \frac{1}{2} ac (T^{\text{NN}}_{\theta}(t_i,x_i))^4 - r^{\text{NN}}_{\theta}(t_i,x_i,\mu_i)) |^2                  \\
                                                      & + \frac{1}{N_\text{int}} \sum_{i=1}^{N_\text{int}}  |\frac{\eps^2}{c\sqrt{\sigma_0}} \partial_t j^{\text{NN}}_{\theta}(t_i,x_i,\mu_i) + \mu_i \partial_x r^{\text{NN}}_{\theta}(t_i,x_i,\mu_i) + \frac{\sigma}{\sqrt{\sigma_0}}j^{\text{NN}}_{\theta}(t_i,x_i,\mu_i) |^2                                                        \\
                                                      & + \frac{1}{N_\text{int}}\sum_{i=1}^{N_\text{int}} | \frac{1}{c} \partial_t \rho^{\text{NN}}_{\theta}(t_i,x_i) + \frac{1}{\sqrt{\sigma_0}} \left \langle \mu \partial_x j^{\text{NN}}_{\theta} \right \rangle(t_i,x_i) + \frac{1}{2}C_v\frac{\partial T}{\partial t}(t_i,x_i)  |^2   \\
                                                      & +\frac{1}{N_\text{int}}\sum_{i=1}^{N_\text{int}} | \eps^2 C_v\frac{\partial T^{\text{NN}}_{\theta}}{\partial t}(t_i,x_i) - \sigma (2 \rho^{\text{NN}}_{\theta}(t_i, x_i) - a c (T^{\text{NN}}_{\theta}(t_i,x_i))^4 ) |^2 ,
                                                      \\
        \mathcal{R}^{\eps}_{\text{constraint}} & = \frac{1}{N_\text{int}}\sum_{i=1}^{N_\text{int}} |\rho^{\text{NN}}_{\theta}(t_i,x_i) - \left \langle r^{\text{NN}}_{\theta} \right \rangle(t_i,x_i) |^2,                                                                                                                            \\
        \mathcal{R}^{\eps}_{\text{boundary}}   & = \frac{1}{N_\text{bdy}} \sum_{i=1}^{N_\text{bdy}}  |\mathcal{B}(r^{\text{NN}}_{\theta} + \frac{\eps}{\sqrt{\sigma_0}} j^{\text{NN}}_{\theta})(t_i,x_i,\mu_i) |^2, \\
        \mathcal{R}^{\eps}_{\text{initial}}    & = \frac{1}{N_0}\sum_{i=1}^{N_0}  | T^{\text{NN}}_{\theta}(0, x_i) - T_0 (x_i)|^2 + \frac{1}{N_0} \sum_{i=1}^{N_0} | (r^{\text{NN}}_{\theta} + \frac{\eps}{\sqrt{\sigma_0}} j^{\text{NN}}_{\theta})(0, x_i,\mu_i) - I_{0}(x_i,\mu_i) |^2.
    \end{aligned}
\end{equation}

\begin{figure}[ht]
    \centering
    \includegraphics[width=0.75\textwidth]{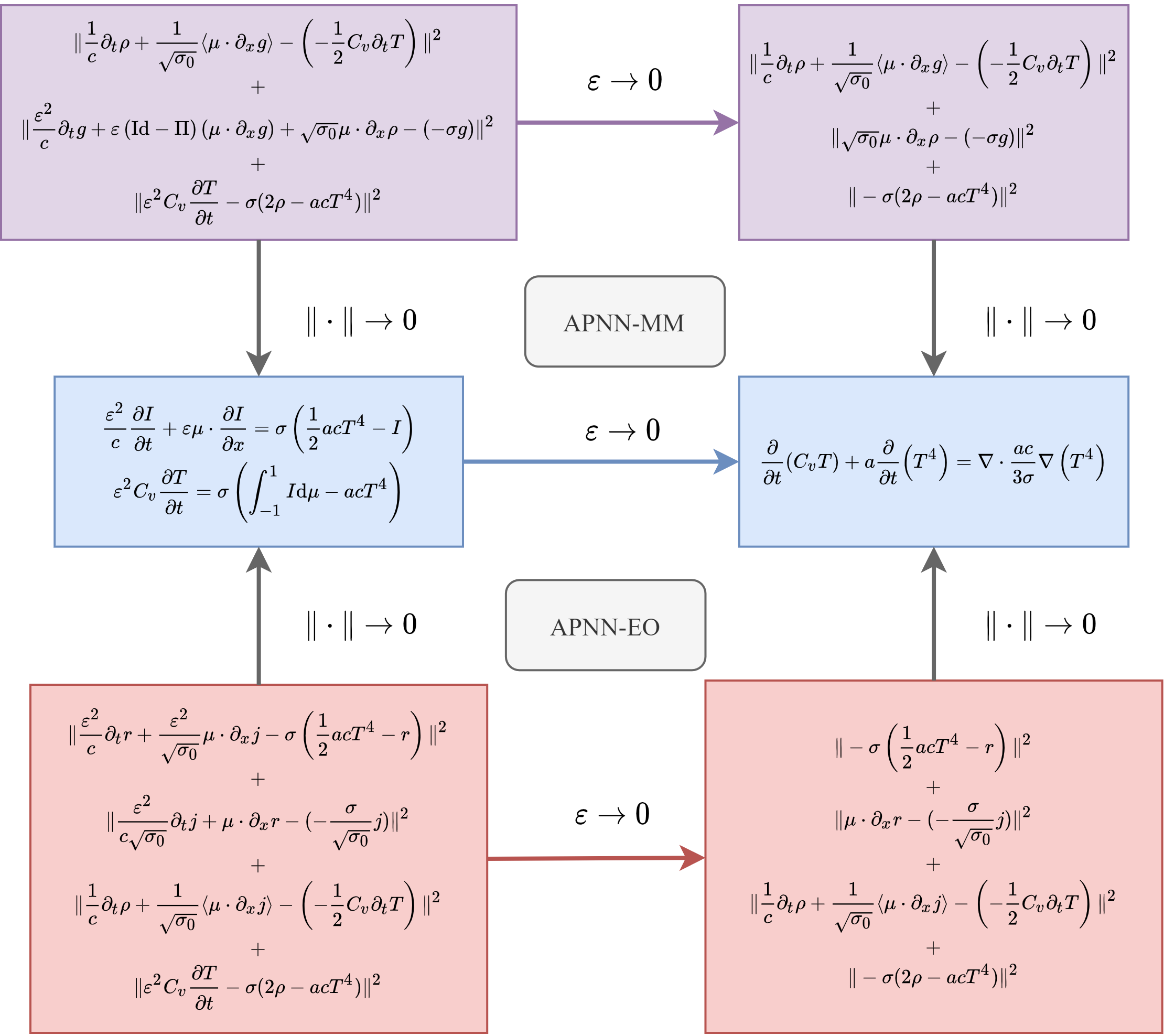}
    \caption{The idea of APNN-MM and APNN-EO for solving the GRTEs.}
    \label{fig:apnn-grte}
\end{figure}

The following diagram~\ref{fig:apnn-grte} illustrates the idea of APNN-MM and APNN-EO methods for solving the GRTEs.

\begin{figure}[ht]
    \centering
    \includegraphics[width=0.8\textwidth]{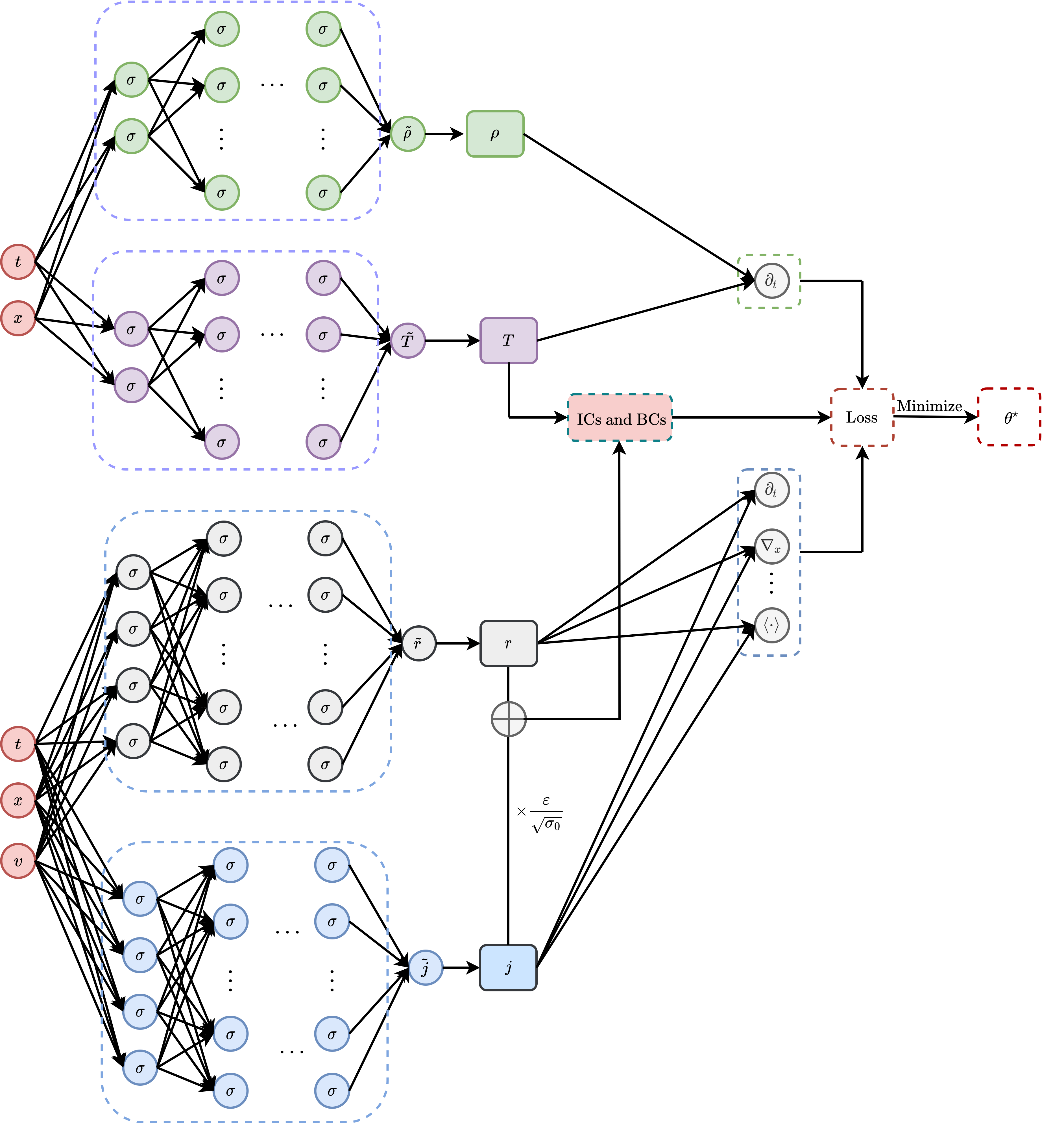}
    \caption{Schematic plot of APNNs based on even-odd decomposition for solving the GRTEs.}
    \label{fig:APNNs}
\end{figure}

Finally, we present a schematic plot depicting our Asymptotic-Preserving Neural Network based on even-odd decomposition for GRTEs in Figure~\ref{fig:APNNs}.

\section{Numerical results}
In this section, we have conducted several numerical experiments to validate the performance of our proposed APNN method for solving GRTEs. As the operators involving $\mu$ in the loss function of APNNs are represented as integrals, we approximate these integrals using the Gauss-Legendre quadrature rule with the number of quadrature points set to 16.

The reference solutions are obtained by spherical harmonics method~\cite{evans1998spherical, kourganov1963basic} and
we will check the relative $\ell^2$ error of the solution $s(x)$ of APNN method, e.g.\ for $1$d case,

\begin{equation}
    \text{error} := \sqrt{
\sum_j |s_{\text{nn}}^{j} - s_{\text{ref}}^j|^2 /
    {\sum_j |s_{\text{ref}}^j|^2}},
\end{equation}
where $s_{\text{nn}}$ is the neural solution approximation, and $s_{\text{ref}}$ is the reference solution.
{Here, the superscript $j$ refers to the index of collocation points at which both the neural approximation and the reference solution are computed.}

\subsection{Experiment setting}
In the experiments, we utilize ResNet with the activation function $\sigma(x) = \text{gelu}(x)$ for all the test problems. 
For the linear transport problem, we employ 2-block ResNets with units of $128$ for $\rho$ and $256$ for $g, r, j$, and $f$. 
{In \cite{wuAPNNv2}, we find that the numerical performance of this structure for the nonlinear problems is better than general MLP.} 
On the other hand, for the stationary and time-dependent nonlinear GRTE problems, we use 3-block ResNets with units of $96$ for $\rho$ and $T$, and $128$ for $g, r, j$.
The spatial domain of interest, denoted as $D := [x_L, x_R]$, covers the interval $[0, 1]$ for the linear transport equation and stationary nonlinear GRTE, while it spans $[0, 2]$ for the time-dependent nonlinear GRTEs.
To train the networks, we utilize the Adam optimizer~\cite{kingma2014adam} version of the gradient descent method with Xavier initialization. 
In the process, hyperparameters, such as neural network architecture, learning rate, and batch size, need to be tuned to achieve a satisfactory level of accuracy~\cite{XDE}.
In each iteration, we use $2048$ random sample points for the domain, $512$ for the boundary, and $1024$ for the initial condition. 
Additionally, to enhance numerical performance, we employ a decreasing annealing schedule for the learning rate.
We use an exponential decay strategy for an initial learning rate $\eta_0 = 10^{-3}$ with a decay rate of $\gamma = 0.96$ and a decay step of $p = 500$ iterations:
\begin{equation*}
    \eta_t = \eta_0 \cdot \gamma^{\lfloor \frac{t}{p} \rfloor},
\end{equation*}
here, the variable $t$ represents the current $t-$th iteration step, and the symbol $\lfloor \cdot \rfloor$ denotes the floor function.

\subsection{Problem 1: APNNs for solving the linear transport equation}
As a warm-up task, we apply APNNs based on micro-macro and even-odd decomposition to solve the 1D linear radiative transfer under the diffusion regime with $\eps = 10^{-3}$, i.e., 
\begin{equation}\label{eq:lte}
    \left\{
    \begin{aligned}
         & \frac{\eps}{c} \partial_t I+ \mu \partial_x I =\frac{\sigma}{\eps} \left(\left \langle I \right \rangle -I\right), \; (x,\mu) \in D \times [-1,1],               \\
         & I(t, x_L, \mu>0) = 1, \; I(t, x_R, \mu<0) = 0, \\
         & I(0, x, v) = 0,
    \end{aligned} 
    \right.
\end{equation}
where $\sigma = 1$. For the corresponding AP loss of the linear transport equation, please refer to~\cite{wuAPNN,wuAPNNv2}.

\begin{figure}[ht]
    \subfigure[{\it Left: PINN}.]
    {
    \includegraphics[width=0.3\textwidth]{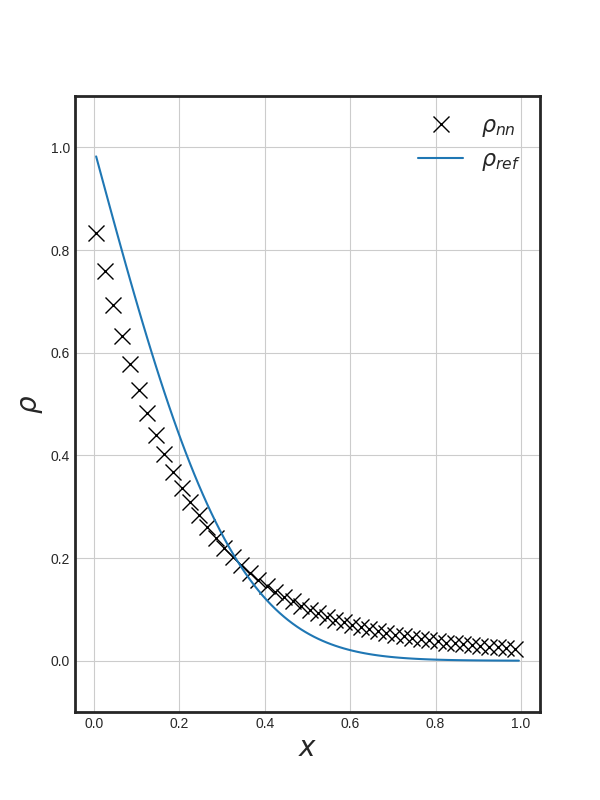}
    }
    \subfigure[{\it Middle: APNN-MM}.]
    {
    \includegraphics[width=0.3\textwidth]{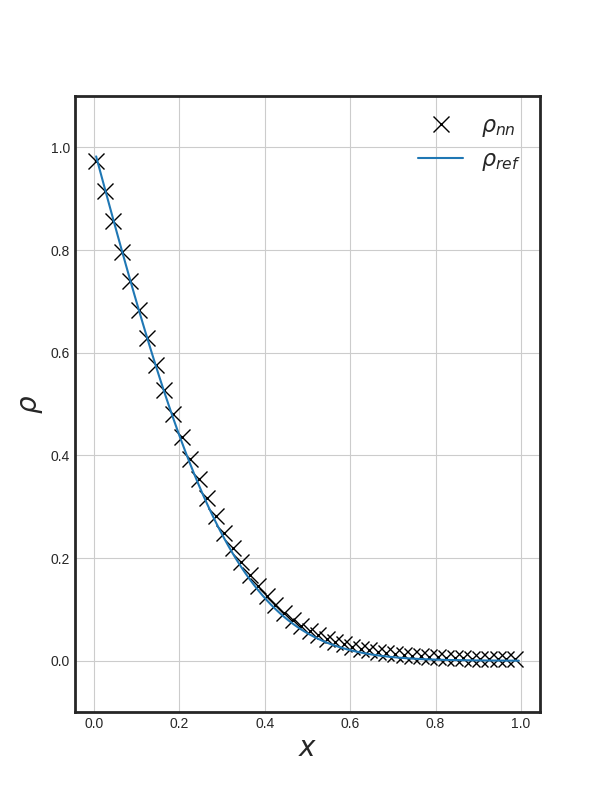}
    }
    \subfigure[{\it Right: APNN-EO}.]
    {
    \includegraphics[width=0.3\textwidth]{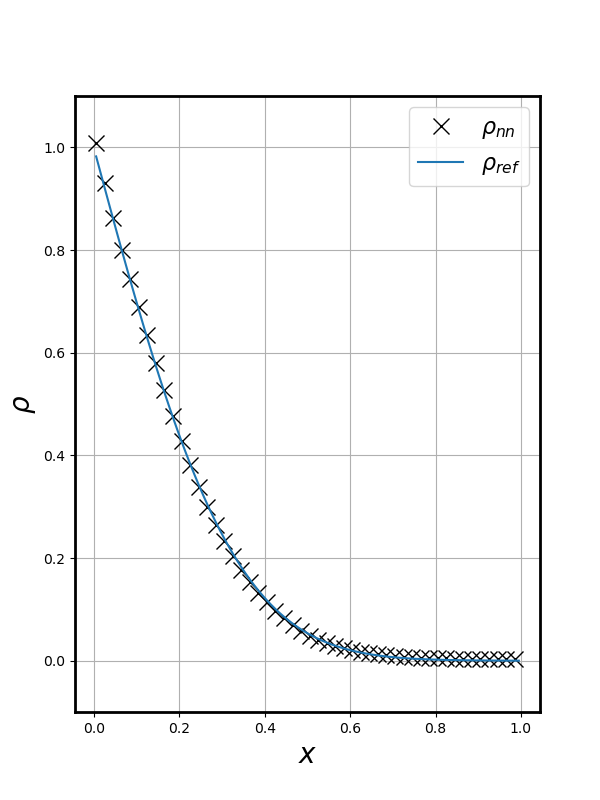}
    }
    
    \caption{\@ Plot of density at time $t = 0.1$ ($\eps = 10^{-3}$): Approximated by PINN, APNNs based on micro-macro and even-odd decomposition (marker) vs. reference solution (line).
    }\label{fig:lte-1e-3}
\end{figure}

Figure~\ref{fig:lte-1e-3} shows the estimated density $\rho$ using DNNs (PINN, APNNs based on micro-macro and even-odd decomposition) in comparison to the reference solution at time $t = 0.1$. 
The relative $\ell^2$ error of PINN, APNNs based on micro-macro and even-odd decomposition are $1.98 \times 10^{-1}, 1.52 \times 10^{-2}, 9.90 \times 10^{-3}$.
The results demonstrate that the approximated solutions obtained through the APNNs exhibit superior accuracy compared to the poor performance of the PINN method.

\subsection{Problem 2: APNNs for solving the stationary nonlinear GRTE}
Next, we focus our attention on solving the steady-state equation of GRTE.
Consider the 1D steady nonlinear gray radiative transfer equations~\cite{lu2022solving,li2022model} given by 
\begin{equation}\label{eq:stationary-grte}
    \left\{
    \begin{aligned}
         & \eps \mu \frac{\partial I}{\partial x} =\sigma \left( acT^4 -I\right), \; (x,\mu) \in D \times [-1,1],               \\
         & \eps^2 \frac{\partial^2 T}{\partial x^2} = \sigma \left(acT^4 - \left \langle I \right \rangle \right), \; x \in D, \\
         & I(x_L, \mu>0) = 1, \; I(x_R, \mu<0) = 0, \\
         & T(x_L) = 1, \; T(x_R) = 0,
    \end{aligned} 
    \right.
\end{equation}
where $a = c = \sigma = 1$. 

The micro-macro system for the stationary GRTEs is written as follows,
\begin{equation}\label{eq:micro-macro-stationary}
    \left\{
    \begin{aligned}
         & \frac{1}{\sqrt{\sigma_0}} \cdot \left \langle \mu \partial_x g \right \rangle = \frac{\partial^2 T}{\partial x^2}, \\
         & \sqrt{\sigma_0} \cdot \mu \partial_x \rho + \eps \left ( \mu \partial_x g - \left \langle \mu \partial_x g \right \rangle \right ) + \sigma g = 0, \\
         & \eps^2 \frac{\partial^2 T}{\partial x^2} = \sigma \left(acT^4 - \rho \right).
    \end{aligned}
    \right.
\end{equation}

\begin{figure}[ht]
    \centering
    \includegraphics[width=0.4\textwidth]{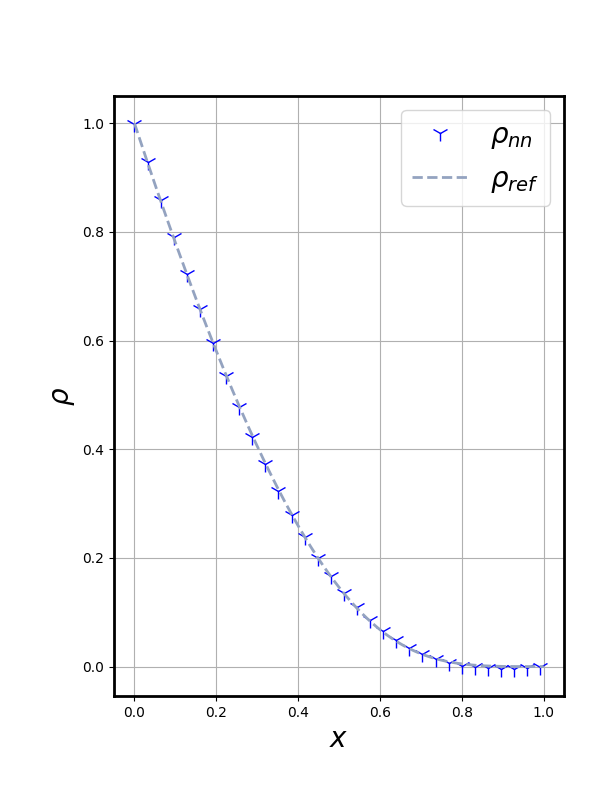}
    \includegraphics[width=0.4\textwidth]{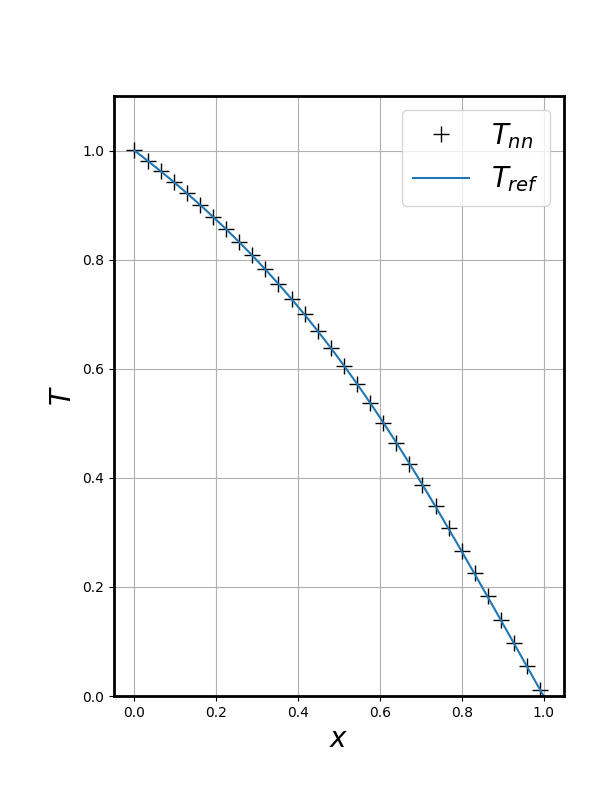}
    \caption{\@ Plot of density and temperature ($\eps = 10^{-3}$): Approximated by APNN based on micro-macro decomposition (marker) vs. reference solutions (line). The relative $\ell^2$ error of $\rho$ and $T$ by APNNs based on micro-macro decomposition are $9.42 \times 10^{-4}, 5.09 \times 10^{-3}$.}
    \label{fig:stationary-micro-macro}
\end{figure}

\begin{figure}[ht]
    \centering
    \includegraphics[width=0.4\textwidth]{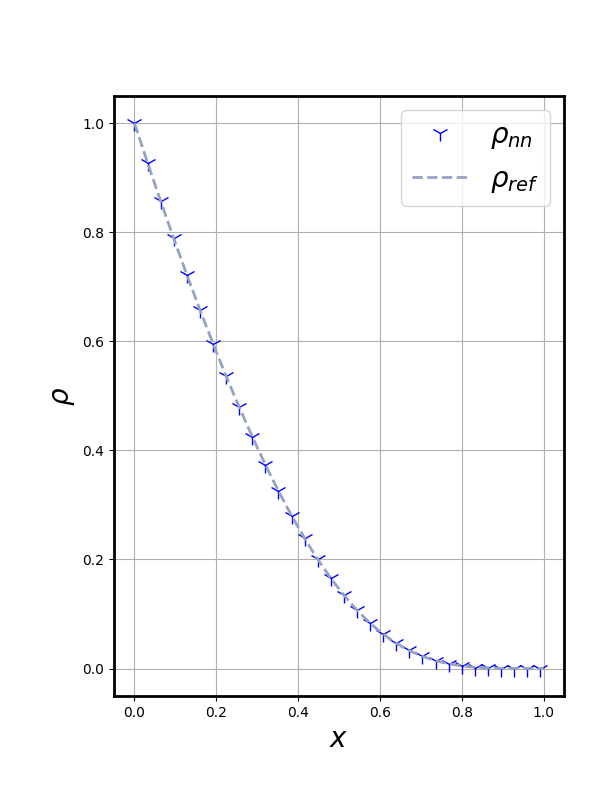}
    \includegraphics[width=0.4\textwidth]{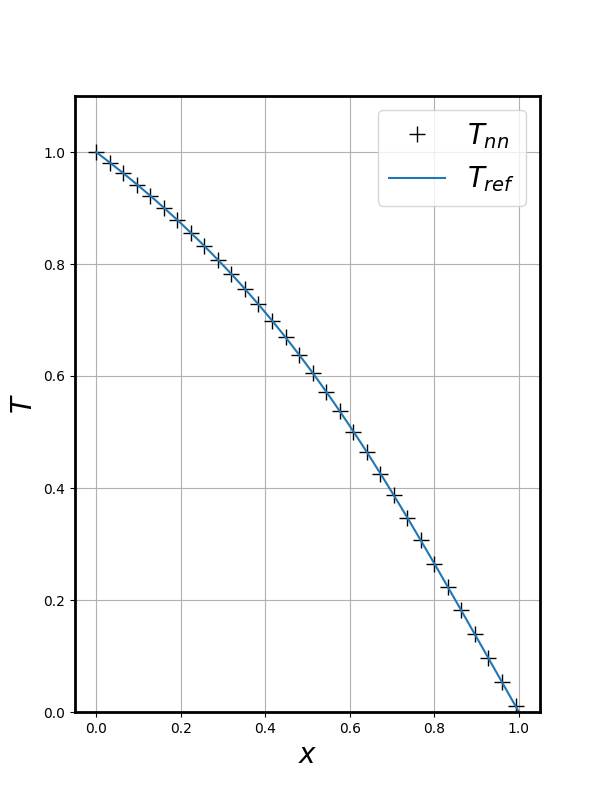}
    \caption{\@ Plot of density and temperature ($\eps = 10^{-3}$): Approximated by APNN based on even-odd decomposition (marker) vs. reference solutions (line). The relative $\ell^2$ error of $\rho$ and $T$ by APNNs based on even-odd decomposition are $6.26 \times 10^{-4}, 9.12 \times 10^{-4}$.}
    \label{fig:stationary-odd-even}
\end{figure}

Furthermore, the stationary GRTEs can be reformulated into the following even-odd system:
\begin{equation}\label{eq:even-odd-stationary}
    \left\{
    \begin{aligned}
         & \frac{\eps^2}{\sqrt{\sigma_0}} \cdot \mu \partial_x j = \sigma \left ( a c T^4 - r \right ), \\
         & \mu \partial_x r = -\frac{\sigma}{\sqrt{\sigma_0}} j, \\
         & \frac{1}{\sqrt{\sigma_0}} \cdot \left \langle  \mu \partial_x j \right \rangle = \frac{\partial^2 T}{\partial x^2}, \\
         & \eps^2 \frac{\partial^2 T}{\partial x^2}  = \sigma \left( acT^4 - \rho \right), \\
         & \rho = \left \langle r \right \rangle.
    \end{aligned}
    \right.
\end{equation}
Figure~\ref{fig:stationary-micro-macro} and \ref{fig:stationary-odd-even} depict the performance achieved by the Asymptotic-Preserving Neural Networks based on micro-macro decomposition and even-odd decomposition, respectively. The relative $\ell^2$ errors of $\rho$ and $T$ obtained through APNNs based on micro-macro decomposition are $9.42 \times 10^{-4}$ and $5.09 \times 10^{-3}$, while for APNNs based on even-odd decomposition, the errors are $6.26 \times 10^{-4}$ and $9.12 \times 10^{-4}$.

The results clearly indicate that the approximate solutions obtained through the Asymptotic-Preserving Neural Networks based on even-odd decomposition outperform those from micro-macro decomposition. One evident reason for this comparison is that the number of parameters in the even-odd decomposition method is higher than the second method.

\subsection{Problem 3: APNNs for solving the time-dependent nonlinear GRTE}
Note that the previous examples solely considered inflow boundary conditions. In this study, we explore the GRTEs with smooth initial conditions transitioning from the kinetic regime ($\eps = 1$) to the diffusion regime ($\eps = 10^{-3}$):
\begin{equation*}
    \left\{
    \begin{aligned}
         & \frac{\eps^2}{c} \frac{\partial I}{\partial t}+\eps \mu \frac{\partial I}{\partial x} =\sigma \left(\frac{1}{2}acT^4 -I\right), \; (t,x,\mu) \in \tau \times D \times [-1,1],               \\
         & \eps^2 C_v\frac{\partial T}{\partial t} = \sigma \left(2 \left \langle I \right \rangle -acT^4 \right), \; (t,x) \in \tau \times D, \\
         & I(t, x_L, \mu) = I(t, x_R, \mu), \\
         & I(0, x, \mu) = \frac{1}{2} ac T(0, x)^4, \; T(0, x) = \frac{3 + \sin(\pi x)}{4},
    \end{aligned}
    \right.
\end{equation*}
where $a = c = 1, C_v = 0.1, \sigma = 10$. The computational region is defined as $[0, 2]$, and periodic boundary conditions are enforced at both ends. The time interval considered in the study is $[0, 0.5]$.

For comparison, we are interested in examining two key quantities: the material temperature denoted as $T_e = T$ and the radiation temperature represented as $T_r = \left ( \frac{1}{ac} \int_{-1}^{1} I \, d\mu \right )^{\frac{1}{4}}$.

\begin{figure}[ht]
    \centering
    \includegraphics[width=0.4\textwidth]{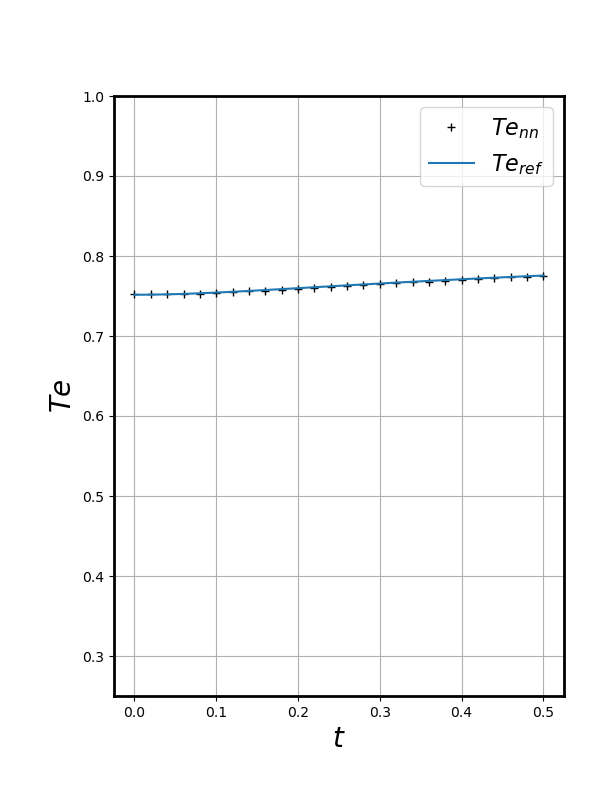}
    \includegraphics[width=0.4\textwidth]{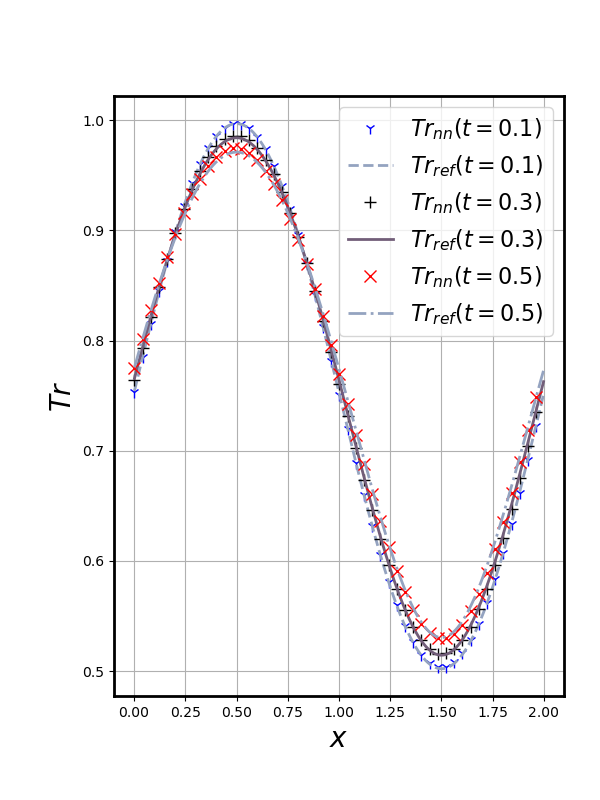}
    \caption{\@ Plot of the approximated material temperature $T_e = T$ at space $x = 0.0025$ and the radiation temperature $T_r$ at time $t = 0.1, 0.3, 0.5$ with APNN based on micro-macro decomposition under the kinetic regime.}
    \label{fig:micro-macro-1e0}
\end{figure}

\begin{figure}[ht]
    \centering
    \includegraphics[width=0.4\textwidth]{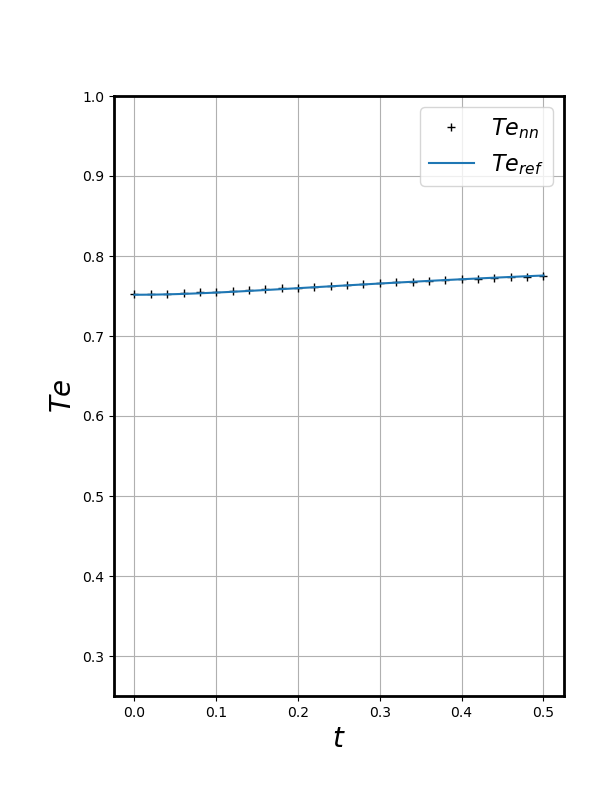}
    \includegraphics[width=0.4\textwidth]{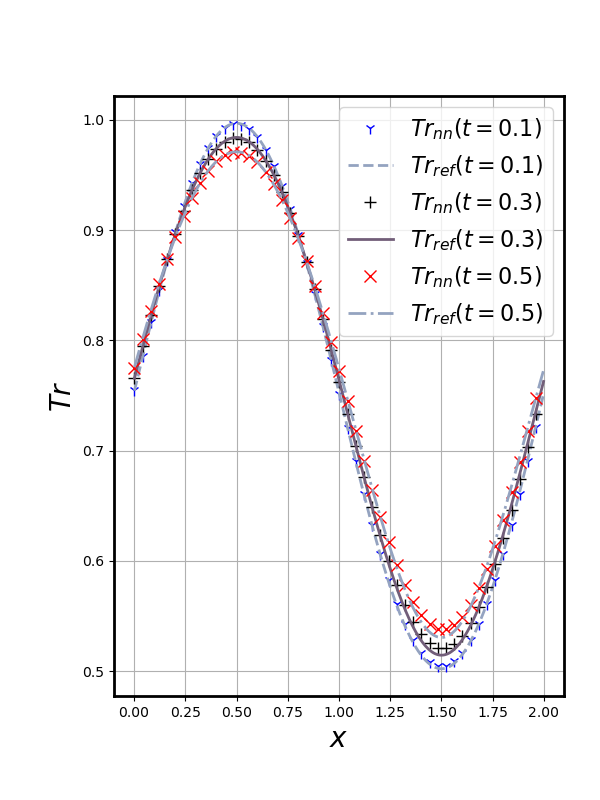}
    \caption{\@ Plot of the approximated material temperature $T_e = T$ at space $x = 0.0025$ and the radiation temperature $T_r$ at time $t = 0.1, 0.3, 0.5$ with APNN based on even-odd decomposition under the kinetic regime.}
    \label{fig:odd-even-1e0}
\end{figure}

Figure~\ref{fig:micro-macro-1e0} and~\ref{fig:odd-even-1e0} illustrates the approximated material temperature $T_e = T$ at space $x = 0.0025$ and the radiation temperature $T_r$ at times $t = 0.1, 0.3, 0.5$ using APNN based on micro-macro decomposition and even-odd decomposition under the kinetic regime.

\begin{figure}[ht]
    \centering
    \includegraphics[width=0.4\textwidth]{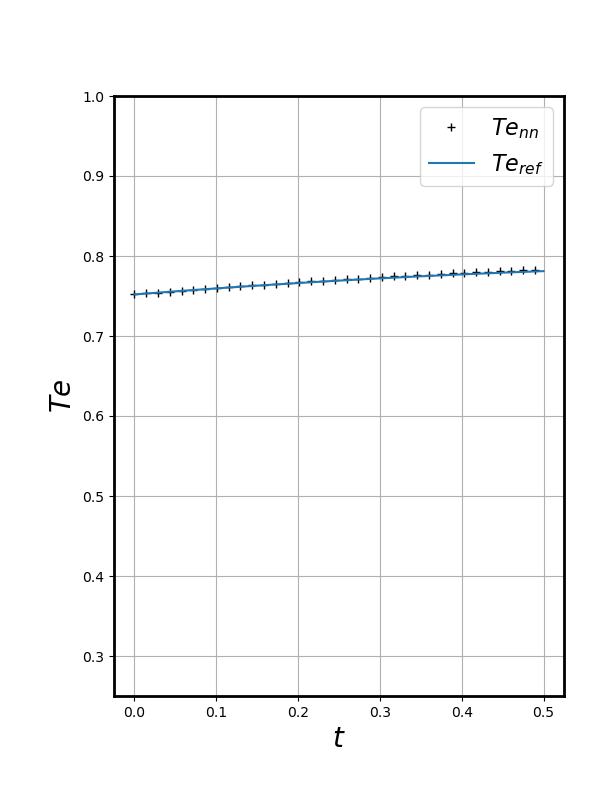}
    \includegraphics[width=0.4\textwidth]{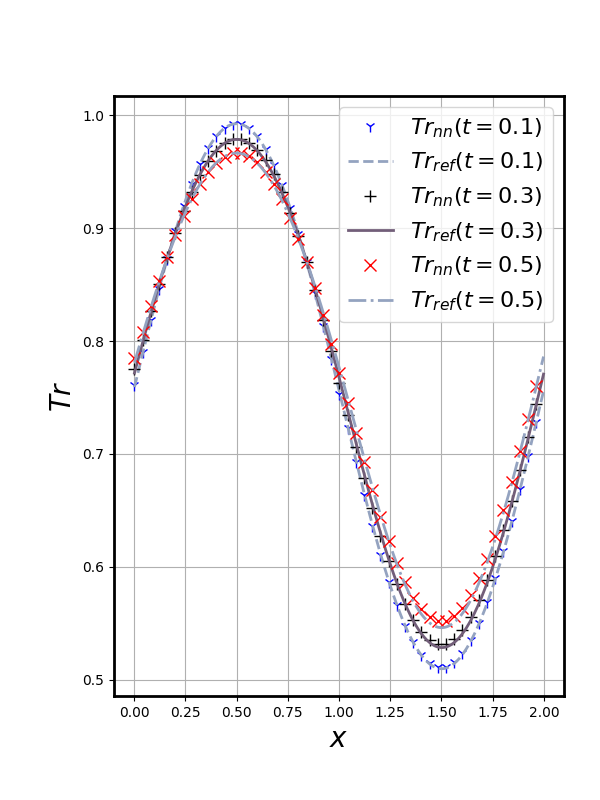}
    \caption{\@ Plot of the approximated material temperature $T_e = T$ at space $x = 0.0025$ and the radiation temperature $T_r$ at time $t = 0.1, 0.3, 0.5$ with APNN based on micro-macro decomposition under the diffusion regime.}
    \label{fig:micro-macro-1e-3}
\end{figure}

\begin{figure}[ht]
    \centering
    \includegraphics[width=0.4\textwidth]{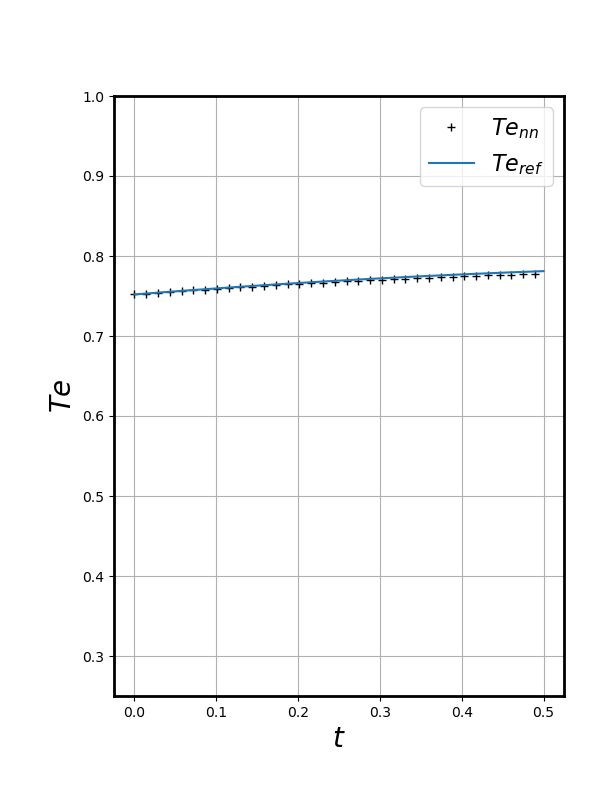}
    \includegraphics[width=0.4\textwidth]{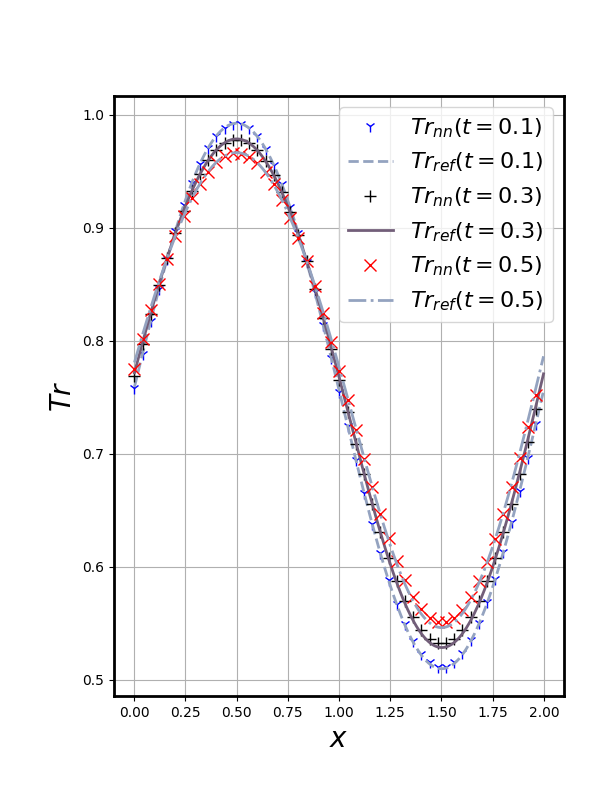}
    \caption{\@ Plot of the approximated material temperature $T_e = T$ at space $x = 0.0025$ and the radiation temperature $T_r$ at time $t = 0.1, 0.3, 0.5$ with APNN based on even-odd decomposition under the diffusion regime.}
    \label{fig:odd-even-1e-3}
\end{figure}

Figure~\ref{fig:micro-macro-1e-3} and~\ref{fig:odd-even-1e-3}  illustrates the inferred material temperature $T_e = T$ at spatial position $x = 0.0025$ and the radiation temperature $T_r$ at time instances $t = 0.1, 0.3, 0.5$, employing the APNN methodology grounded on the micro-macro decomposition and even-odd decomposition approach, within the diffusion regime.

 Table~\ref{tab:error-micro-macro} documents the relative $\ell^2$ error of $T_e$ and $T_r$ by APNN based on micro-macro decomposition at time $t = 0.1, 0.2, 0.3, 0.4, 0.5$ for the kinetic regime ($\eps = 1$) and diffusion regime ($\eps = 10^{-3}$).

\begin{small}
    \begin{table}[htbp]
        \caption{The relative $\ell^2$ error of $T_e$ and $T_r$ by APNN based on micro-macro decomposition at time $t = 0.1, 0.2, 0.3, 0.4, 0.5$ for the kinetic regime ($\eps = 1$) and diffusion regime ($\eps = 10^{-3}$).}
        \label{tab:error-micro-macro}
        \centering
        \begin{tabular}{ccccccccc}
            \toprule[1pt]
            \noalign{\smallskip}
            \multirow{2}*{\diagbox{{$\eps$}}{Error}}
             & \multicolumn{6}{c}{APNN-MM}  \\
             & \multicolumn{1}{c}{$T_e$} & \multicolumn{1}{c}{$T_r (t = 0.1)$} & \multicolumn{1}{c}{$T_r (t = 0.2)$} & \multicolumn{1}{c}{$T_r (t = 0.3)$} & \multicolumn{1}{c}{$T_r (t = 0.4)$} & \multicolumn{1}{c}{$T_r (t = 0.5)$}  \\
            \noalign{\smallskip}
            \midrule[1pt]
            \noalign{\smallskip}
            \multirow{1}*{$1$}
             & $9.25{\text{e-}4}$ & $7.01{\text{e-}4}$ & $1.50{\text{e-}3}$ & $2.43{\text{e-}3}$ & $3.42{\text{e-}3}$ & $4.37{\text{e-}3}$ \\
            \multirow{1}*{{{$0.001$}}}
             & $3.11{\text{e-}3}$ & $1.22{\text{e-}3}$ & $2.14{\text{e-}3}$ & $2.81{\text{e-}3}$ & $3.26{\text{e-}3}$ & $3.57{\text{e-}3}$                   \\
            \noalign{\smallskip}
            \bottomrule[1pt]
        \end{tabular}
    \end{table}
\end{small}

Table~\ref{tab:error-odd-even} documents the relative $\ell^2$ error of $T_e$ and $T_r$ by APNN based on even-odd decomposition at time $t = 0.1, 0.2, 0.3, 0.4, 0.5$ for the kinetic regime ($\eps = 1$) and diffusion regime ($\eps = 10^{-3}$).

\begin{small}
    \begin{table}[htbp]
        \caption{The relative $\ell^2$ error of $T_e$ and $T_r$ by APNN based on even-odd decomposition at time $t = 0.1, 0.2, 0.3, 0.4, 0.5$ for the kinetic regime ($\eps = 1$) and diffusion regime ($\eps = 10^{-3}$).}
        \label{tab:error-odd-even}
        \centering
        \begin{tabular}{ccccccccc}
            \toprule[1pt]
            \noalign{\smallskip}
            \multirow{2}*{\diagbox{{$\eps$}}{Error}}
             & \multicolumn{6}{c}{APNN-EO}  \\
             & \multicolumn{1}{c}{$T_e$} & \multicolumn{1}{c}{$T_r (t = 0.1)$} & \multicolumn{1}{c}{$T_r (t = 0.2)$} & \multicolumn{1}{c}{$T_r (t = 0.3)$} & \multicolumn{1}{c}{$T_r (t = 0.4)$} & \multicolumn{1}{c}{$T_r (t = 0.5)$}  \\
            \noalign{\smallskip}
            \midrule[1pt]
            \noalign{\smallskip}
            \multirow{1}*{$1$}
             & $8.70{\text{e-}4}$ & $1.47{\text{e-}3}$ & $2.65{\text{e-}3}$ & $3.44{\text{e-}3}$ & $3.77{\text{e-}3}$ & $3.78{\text{e-}3}$ \\
            \multirow{1}*{$0.001$}
             & $3.80{\text{e-}3}$ & $8.83{\text{e-}4}$ & $1.61{\text{e-}3}$ & $2.28{\text{e-}3}$ & $2.90{\text{e-}3}$ & $3.46{\text{e-}3}$  \\
            \noalign{\smallskip}
            \bottomrule[1pt]
        \end{tabular}
    \end{table}
\end{small}

\subsection{Problem 4: APNNs for solving another time-dependent GRTE } 
we have solved an initial boundary value incompatibility problem, where the opacity $(\sigma=10 \mathrm{cm}^{-1} )$ is temperature-independent, and the heat capacity $( C_{v} = 1 \; \mathrm{GJ} / \mathrm{cm}^3 / \mathrm{KeV} )$ in the context of one-dimensional time-dependent GRTEs. This slab, with a thickness of 0.25 $\mathrm{cm}$, starts in an equilibrium state at 1 $ \mathrm{KeV}$ and is subject to reflection and incident Planckian source conditions at the left and right boundaries, respectively.
\begin{equation*}
    \left\{
    \begin{aligned}
        &\frac{\eps^2}{c} \frac{\partial I}{\partial t}+\eps \mu \frac{\partial I}{\partial x}=\sigma\left(\frac{1}{2} a c T^4-I\right) ,\; (t, x, \mu) \in \tau \times D \times[-1,1], \\ 
        &\eps^2 C_v \frac{\partial T}{\partial t}=\sigma\left(2\langle I\rangle-a c T^4\right) ,\; (t, x) \in \tau \times D , \\
        &I(t, 0, \mu>0)=I(t, 0,-\mu), \\ 
        &I(t, 0.25, \mu<0)=\frac{1}{2} a c(0.1)^4 , \\ 
        &I(0, x, \mu)=\frac{1}{2} a c T(0, x)^4, \quad T(0, x) = 1 ,
    \end{aligned}
    \right.
\end{equation*}
where $\eps=1, a=0.01372, c=29.98, C_v=1, \sigma=10$. The computational region is defined as [0,0.25], and the time interval considered in the study is $[0,1]$.

\begin{figure}[ht]
    \centering
    \includegraphics[width=0.4\textwidth]{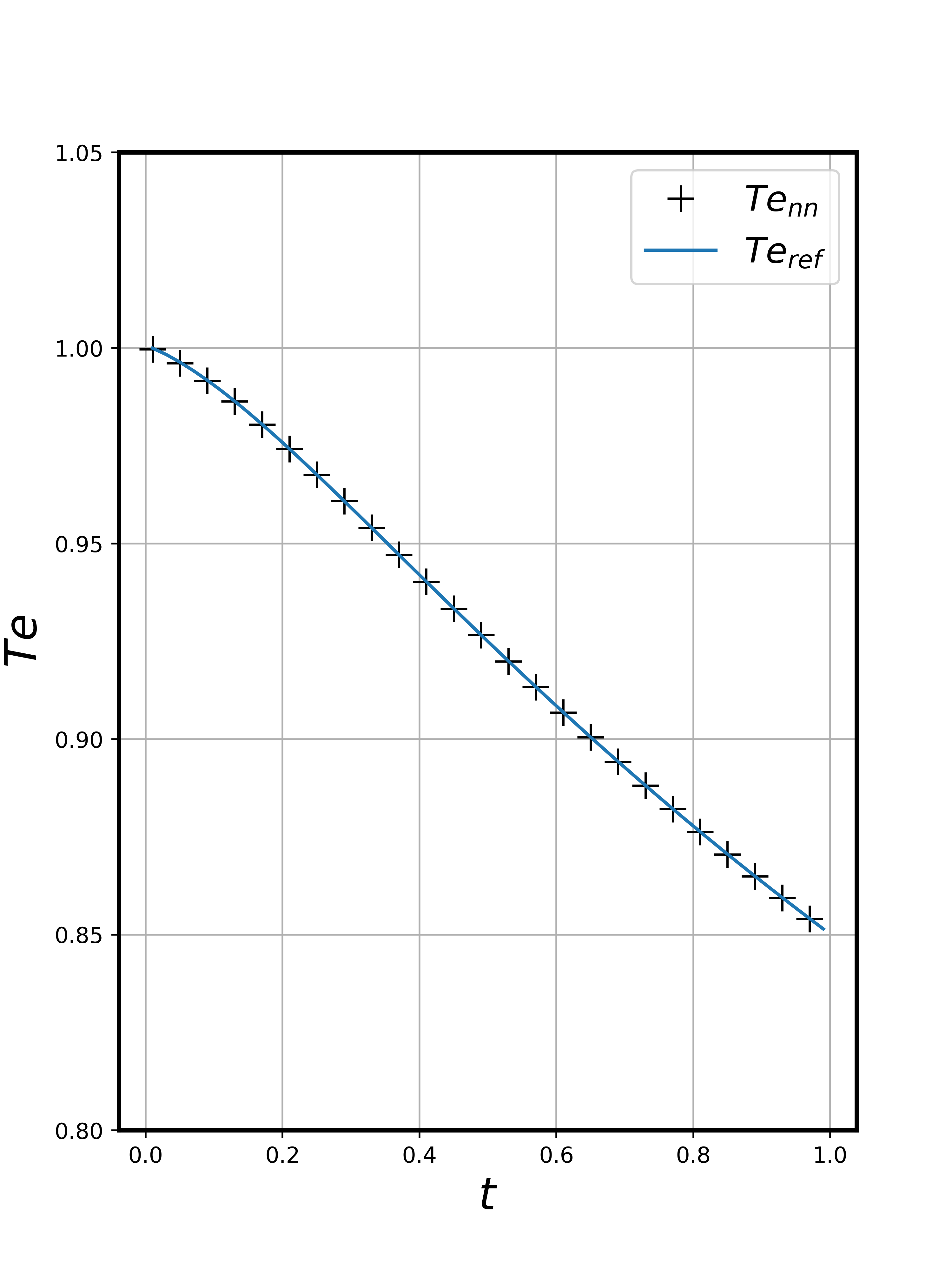}
    \includegraphics[width=0.4\textwidth]{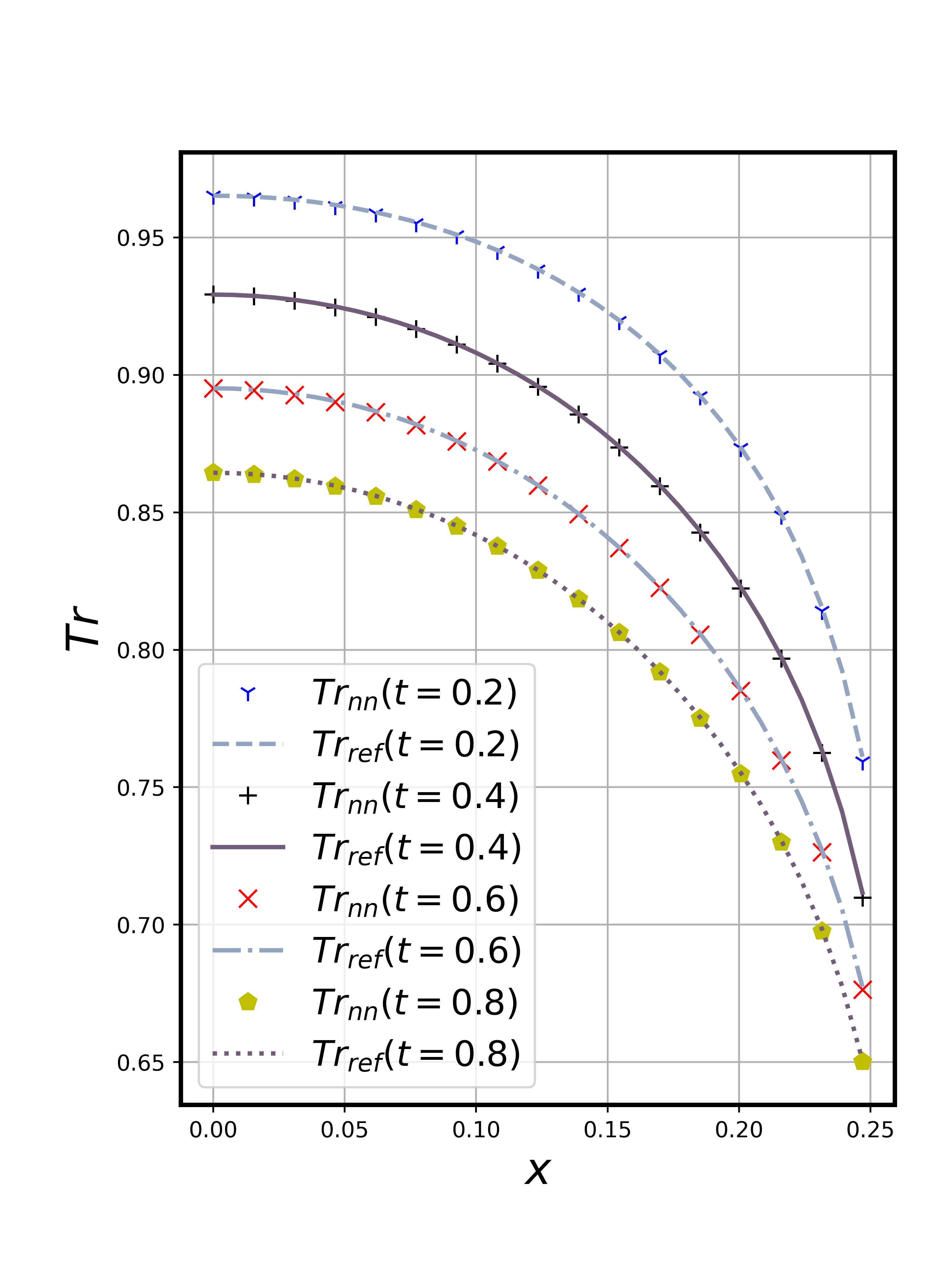}
    \caption{\@ Plot of the approximated material temperature $T_e = T$ at space $x = 0.0025$ and the radiation temperature $T_r$ at time $t = 0.2, 0.4, 0.6, 0.8$ with APNN based on even-odd decomposition under the kinetic regime.}\label{fig:odd_even_1e0}
\end{figure}

Figure~\ref{fig:odd_even_1e0} illustrates the estimation of the approximate material temperature $T_e=T$ and radiation temperature $T_r$ using the APNN method based on even-odd decomposition in a dynamic state. The estimations are presented at spatial location $x=0.0025$ and time instances $t=0.2, 0.4, 0.6, 0.8$.

Table~\ref{tab:error2-odd-even} records the relative $\ell^2$ errors of $T_r$ at time instances $t=0.2, 0.4, 0.6, 0.8$, and of $T_e$ at spatial location $x=0.0025$ for the APNN method based on even-odd decomposition in the dynamic state ($\eps=1$).

\begin{figure}[ht]
    \centering
    \includegraphics[width=0.8\textwidth]{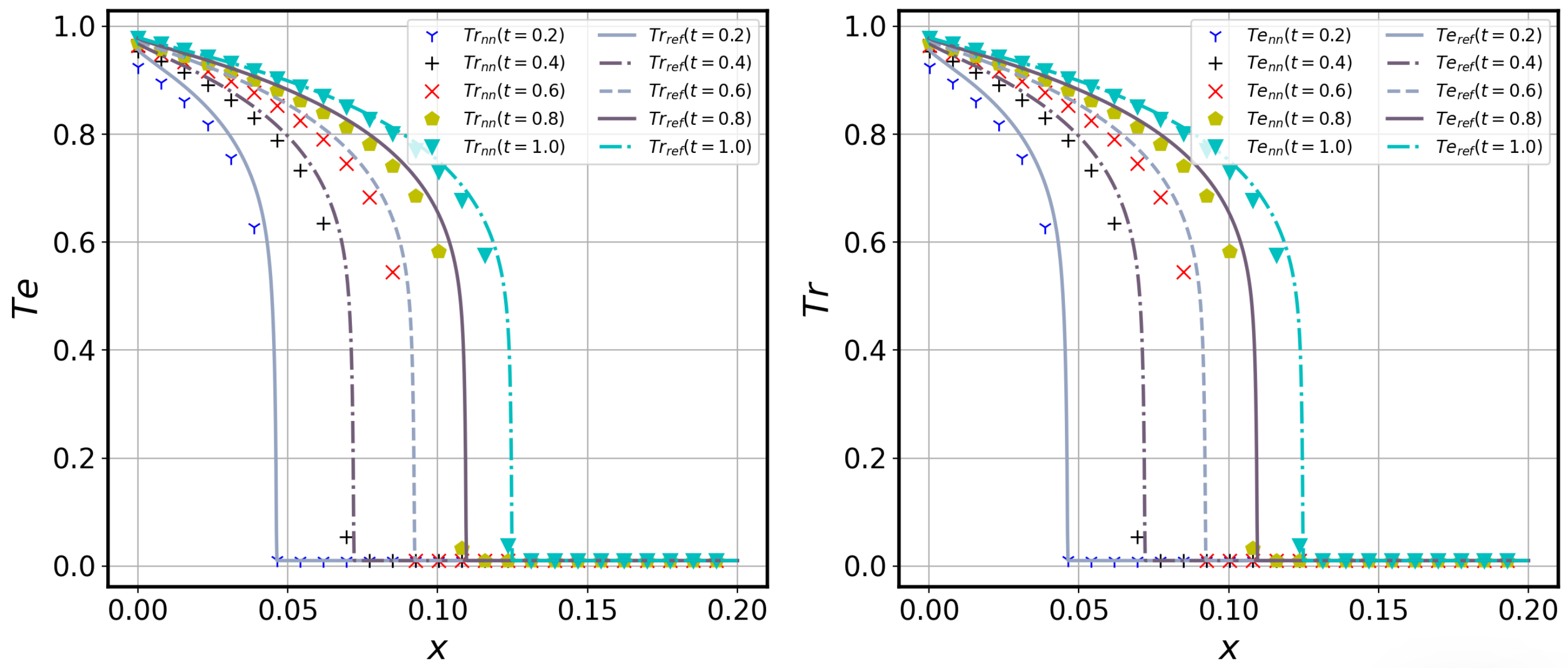}
    \caption{\@ Plot of the approximated material temperature $T_e = T$  and the radiation temperature $T_r$ at time $t = 0.2,0.4,0.6,0.8,1.0$ with APNN based on even-odd decomposition.}
    \label{fig:marshak_wave}
\end{figure}

\begin{small}
    \begin{table}[htbp]
        \caption{The relative $\ell^2$ error of $T_e$ and $T_r$ by APNN based on even-odd decomposition at time $t =0.2, 0.4, 0.6, 0.8$ for the kinetic regime ($\eps = 1$).}
        \label{tab:error2-odd-even}
        \centering
        \begin{tabular}{ccccccccc}
            \toprule[1pt]
            \noalign{\smallskip}
            \multirow{2}*{\diagbox{{$\eps$}}{Error}}
             & \multicolumn{5}{c}{APNN-EO}  \\
             & \multicolumn{1}{c}{$T_e$}  & \multicolumn{1}{c}{$T_r (t = 0.2)$} & \multicolumn{1}{c}{$T_r (t = 0.4)$} & \multicolumn{1}{c}{$T_r (t = 0.6)$} & \multicolumn{1}{c}{$T_r (t = 0.8)$}  \\
            \noalign{\smallskip}
            \midrule[1pt]
            \noalign{\smallskip}
            \multirow{1}*{$1$}
             & $1.05{\text{e-}4}$ & $6.54{\text{e-}4}$ & $6.39{\text{e-}4}$ & $5.00{\text{e-}4}$ & $4.81{\text{e-}4}$  \\
            
            \noalign{\smallskip}
            \bottomrule[1pt]
        \end{tabular}
    \end{table}
\end{small}
{
\subsection{Problem 5: APNNs for solving Marshak wave problem }
Finally, we study a problem with a temperature-dependent absorption coefficient, namely the Marshak wave problem: 
\begin{equation*}
\left\{
\begin{aligned}
&\frac{\eps^2}{c} \frac{\partial I}{\partial t}+\varepsilon \mu \frac{\partial I}{\partial x}=\sigma\left(\frac{1}{2} a c T^4-I\right),(t, x, \mu) \in \tau \times D \times[-1,1], \\
&\eps^2 C_v \frac{\partial T}{\partial t}=\sigma\left(2\langle I\rangle-a c T^4\right),(t, x) \in \tau \times D, \\
&I\left(t, 0, \mu>0\right)=\frac{1}{2} a c , \\
&I\left(t, 0.2, \mu<0\right)=\frac{1}{2} a c (0.01)^4, \\
&I(0, x, \mu)=\frac{1}{2} a c (0.01)^4 ,
\end{aligned}\right.
\end{equation*}
where $c=29.98, a=0.01372$, and $C_v=0.3$. The absorption coefficient is given by $\sigma(T)=\frac{30}{T^3}$. Inflow boundary conditions are applied at both ends. Due to the rapid variation in material temperature at the initial time, the problem is solved over the time interval $\tau=[0.1 \mathrm{~ns}, 1 \mathrm{~ns}]$.
Figure~\ref{fig:marshak_wave} illustrates the estimation of material temperature $T_e=T$ and radiation temperature $T_r$ using the APNN method based on even-odd decomposition under ( $\eps=1$ ). The estimations are presented at time points $t=0.2,0.4,0.6,0.8,1.0$. Table~\ref{tab:marshak_wave} records the relative $\ell^2$ errors of the radiation temperature $T_r$ and material temperature $T_e$ at these time points. The results show that we have successfully simulated the propagation of the Marshak wave, but there is still a deviation in the prediction of the wave front.}

\begin{small}
    \begin{table}[htbp]
        \caption{The relative $\ell^2$ error of $T_e$ and $T_r$ by APNN based on even-odd decomposition at time $t =0.2, 0.4, 0.6, 0.8,1.0$ for the kinetic regime ($\eps = 1$).}
        \label{tab:marshak_wave}
        \centering
        \begin{tabular}{cccccccc}
            \toprule[1pt]
            \noalign{\smallskip}
            \multirow{2}*{\diagbox{{$T$}}{Error}}
             & \multicolumn{5}{c}{APNN-EO}  \\
             & \multicolumn{1}{c}{$t = 0.2$} & \multicolumn{1}{c}{$t = 0.4$} & \multicolumn{1}{c}{$t = 0.6$} & \multicolumn{1}{c}{$t = 0.8$} & \multicolumn{1}{c}{$t = 1.0$}  \\
            \noalign{\smallskip}
            \midrule[1pt]
            \noalign{\smallskip}
            \multirow{1}*{$T_e$}
             & $1.40{\text{e-}1}$ & $1.35{\text{e-}1}$ & $1.07{\text{e-}1}$ & $8.96{\text{e-}2}$ & $7.68{\text{e-}2}$ \\
            \multirow{1}*{$T_r$}
             & $1.46{\text{e-}1}$ & $1.38{\text{e-}1}$ & $1.09{\text{e-}1}$ & $9.15{\text{e-}2}$ & $7.89{\text{e-}2}$  \\
            \noalign{\smallskip}
            \bottomrule[1pt]
        \end{tabular}
    \end{table}
\end{small}

\section{Conclusion}
In this research article, we present a new
APNN approach for the nonlinear gray radiative transfer equations, which is built upon the even-odd decomposition technique. 
This novel method aims to efficiently solve the nonlinear gray radiative transfer equations. 
Our method introduces an auxiliary deep neural network, distinct from the micro-macro decomposition-based APNN method for GRTEs, and is designed to relax the strict conservation requirements while maintaining uniform stability concerning the small Knudsen number. 
As a result, the neural network solution converges uniformly to {the diffusion limit solution}, ensuring reliable results even at small scales.
To validate the superiority of our proposed method, we conducted several experiments. Initially, we applied the Physics-Informed Neural Network and APNN methods to solve the linear transport equation under the diffusion regime. 
The results revealed that PINN failed to accurately resolve this simplest GRTE with small scales, further highlighting the advantages of our APNN method.
Additionally, we tested the stationary nonlinear GRTE with inflow boundary conditions and a Knudsen number of $\eps = 10^{-3}$ using APNNs based on both micro-macro and even-odd decomposition. 
These two APNN methods demonstrated excellent performance in handling the problem.
Then we tackled the time-dependent nonlinear GRTEs with different boundary conditions, considering both the kinetic regime ($\eps = 1$) and the diffusion regime ($\eps = 10^{-3}$) using APNNs based on micro-macro and even-odd decomposition. 
{Finally, we simulate a relatively difficult Marshak wave problem, which is a problem where the absorption coefficient varies with temperature. We successfully simulate the forward trajectory of the wave. However, there is still a deviation in the prediction of the wave front, and there is significant room for improvement in accuracy.}
The numerical results exhibited the efficacy of our proposed APNN method in solving these complex GRTEs under various conditions.
In conclusion, our research showcases the remarkable capabilities of the novel APNN method, which effectively handles the solution of nonlinear GRTEs and achieves impressive results across different regimes and boundary conditions.

\section*{Acknowledgments}

The work of Keke Wu is supported by the China Postdoctoral Science Foundation under Grant Number 2025M773105 and the Jiangsu Funding Program for Excellent Postdoctoral Talent.
Wengu Chen acknowledges partial support from the National Natural Science Foundation of China (NSFC) under Grant No. 12271050, and the Foundation of National Key Laboratory of Computational Physics (Grant No. 6142A05230503).
Zheng Ma is supported by the National Natural Science Foundation of China (Grant Nos. 12201401 and 92270120), as well as by the Beijing Institute of Applied Physics and Computational Mathematics (Grant No. HX02023-60).

\bibliographystyle{unsrt} 
\bibliography{references}
\end{document}